\documentclass{amsart}

\usepackage{amssymb}

\newtheorem{The}{Theorem}[section]
\newtheorem{Lem}[The]{Lemma}
\newtheorem{Pro}[The]{Proposition}
\newtheorem{Cor}[The]{Corollary}
\newtheorem{Def}[The]{Definition}
\newtheorem{Rem}[The]{Remark}

\newcommand{\N}{\mathbb{N}}
\newcommand{\Z}{\mathbb{Z}}
\newcommand{\R}{\mathbb{R}}
\newcommand{\C}{\mathbb{C}}

\newcommand{\G}{\widehat{G}}
\newcommand{\La}{\Lambda_s^{r_1,r_2}}
\newcommand{\la}{\lambda_s^{r_1,r_2}}
\newcommand{\lan}{\lambda_s^{r_1,r_2}(\cdot,n)}
\newcommand{\Ra}{\mathcal{R}}
\newcommand{\Lap}{\mathrm{L}_{\mathbf{T}}}

\begin{document}

\title[On the Natural Representation of $S(\Omega)$ into
$L^2(\mathcal{P} (\Omega))$]{On the Natural Representation of
$S(\Omega)$ into $L^2(\mathcal{P} (\Omega))$: \\ Discrete
Harmonics and Fourier Transform}

\author[Marco and Parcet]{Jos\'{e} Manuel Marco and Javier Parcet$^{\dag}$}

\address{Department of Mathematics, Universidad Aut\'{o}noma de
Madrid}

\email{javier.parcet@uam.es}

\date{}

\footnote{$^{\dag}$ Partially supported by Project PB97-0030 of
DGES, Spain.} 

\begin{abstract}
Let $\Omega$ denote a non-empty finite set. Let $S(\Omega)$ stand
for the symmetric group on $\Omega$ and let us write $\mathcal{P}
(\Omega)$ for the power set of $\Omega$. Let $\rho: S(\Omega)
\rightarrow U(L^2(\mathcal{P} (\Omega)))$ be the left unitary
representation of $S(\Omega)$ associated with its natural action
on $\mathcal{P} (\Omega)$. We consider the algebra consisting of
those endomorphisms of $L^2(\mathcal{P} (\Omega))$ which commute
with the action of $\rho$. We find an attractive basis
$\mathcal{B}$ for this algebra. We obtain an expression, as a
linear combination of $\mathcal{B}$, for the product of any two
elements of $\mathcal{B}$. We obtain an expression, as a linear
combination of $\mathcal{B}$, for the adjoint of each element of
$\mathcal{B}$. It turns out the Fourier transform on $\mathcal{P}
(\Omega)$ is an element of our algebra; we give the matrix which
represents this transform with respect to $\mathcal{B}$.
\end{abstract}

\maketitle

\section*{Introduction} \label{section1}

Let $\Omega$ be a finite set of $n$ elements. If we denote by $G$
the symmetric group $S(\Omega)$ of permutations of $\Omega$ and by
$X$ the power set $\mathcal{P} (\Omega)$ of $\Omega$, then the
natural action of $G$ on $X$ leads to the associated left
representation $\rho: G \rightarrow U \big( L^2(X) \big)$ given by
$\big( \rho(g) \psi \big) (x) = \psi(g^{-1} (x))$. The aim of this
paper is to study the $\star$--algebra $\mbox{End}_G (L^2(X))$ of
intertwining operators for $\rho$. That is, the algebra of
endomorphisms of the Hilbert space $L^2(X)$ which commute with the
action of $\rho$. The partition of $X$ into orbits of $G$, $X_r =
\{x \in X: |x| = r\}$ $(0 \le r \le n)$, gives rise to a family of
subspaces $L^2(X_r)$ of $L^2(X)$ invariant under the action of
$\rho$. Each $X_r$ is a finite symmetric space with respect to $G$
and there exists a family of inequivalent irreducible
representations $\pi_s: G \rightarrow U(\mathcal{V}_s)$ $(0 \le s
\le [n/2])$, such that the following holds $$L^2(X_r) =
\bigoplus_{s=0}^{r \wedge (n-r)} L^2(X_r)_s.$$ This is a
consequence of the so-called Young rule, see \cite{BI} (p. 212,
theorem 2.5) or \cite{Di} (p. 138-139). Here $\wedge$ stands for
min and $L^2(X_r)_s$ stands for the $G$-invariant subspace of
$L^2(X_r)$ equivalent to $\mathcal{V}_s$. Applying Schur's lemma,
we express the algebra as a direct sum of $1$-dimensional
subspaces as follows. Writing $N(r_1,r_2) = r_1 \wedge (n-r_1)
\wedge r_2 \wedge (n-r_2)$, we have $$\mbox{End}_G (L^2(X)) =
\bigoplus_{0 \le r_1, r_2 \le n} \bigoplus_{s=0}^{N(r_1,r_2)}
\mbox{Hom}_{G} \big( L^2(X_{r_1})_s, L^2(X_{r_2})_s \big).$$

\noindent Taking non zero elements $\La \in \mbox{Hom}_{G} \big(
L^2(X_{r_1})_s, L^2(X_{r_2})_s \big)$ we obtain a basis
$$\mathcal{B} = \big\{ \La: \ \ 0 \le r_1, r_2 \le n, \ \ 0 \le s
\le N(r_1, r_2) \big\}$$ of $\mbox{End}_G (L^2(X))$ which is
orthogonal with respect to the Hilbert-Schmidt inner product.
Since each operator $\La$ commutes with the action of $\rho$, it
is obvious that its kernel $\la(x_2,x_1)$ is constant on
$\theta_k^{r_1,r_2} = \{(x_2,x_1) \in X_{r_2} \times X_{r_1}: |x_2
\setminus x_1| = k\}$, where $0 \vee (r_2-r_1) \le k \le (n-r_1)
\wedge r_2$ and $\vee$ stands for max. We will see that the common
value $\la(k)$ (the evaluation at $k = |x_2 \setminus x_1|$ of the
kernel $\la(x_2,x_1)$) at $\theta_k^{r_1,r_2}$ is given by a Hahn
polynomial. In this paper we provide three expressions for these
polynomials. One of them seems to be new and requires the use of
the Radon transforms, which are operators from $L^2(X_{r_1})$ to
$L^2(X_{r_2})$ that commute with the action of $\rho$. The other
expressions for the Hahn polynomials arise from the use of a
discrete Laplacian operator and the theory of orthogonal
polynomials of hypergeometric type. These last ones can be reduced
to well-known expressions in the range $r_2 \ge r_1$ (see
\cite{NU}) and can be regarded as a symmetryzation of those. These
last two are included in our paper in order to facilitate our
computations. In Theorem \ref{theorem2} we write the products
$\Lambda_s^{r_2,r_3} \circ \La$ in terms of the basis
$\mathcal{B}$. This is the main result of the paper. Its proof
uses the Radon transforms and a characterization of spherical
functions on symmetric spaces which we enunciate at the end of
section \ref{section2}. Finally, using the well-known abelian
group structure on the power set $X$ given by the symmetric
difference operator, we study the associated Fourier transform
$\mathcal{F}_X$ on $X$. We show that it can be considered as a
member of the algebra $\mbox{End}_G (L^2(X))$. Then we apply our
results to this particular case writing the Fourier transform
$\mathcal{F}_X$ as a linear combination of the operators $\La$. We
show that the coefficients of $\mathcal{F}_X$ with respect to
$\mathcal{B}$ can be expressed in terms of the Krawtchouk
polynomials, see theorem \ref{theorem4}. This analysis of
$\mathcal{F}_X$ was one of the motivations of this paper.

The results we present here can be analyzed in terms of
distance-regular graphs, see \cite{BI} or \cite{DL}. Given a
connected distance-regular graph $Y$, with distance $\partial$,
there are two associated algebras. First, the Bose-Mesner algebra
of operators on $L^2(Y)$ whose kernel is a function of
$\partial(x,y)$. Second, the Terwilliger algebra, which is defined
by fixing a point $x_0 \in Y$ and taking the algebra generated by
the operators on $L^2(Y)$ whose kernel is a function of $(\partial
(x,x_0), \partial (x,y), \partial (y,x_0))$. The reader is
referred to \cite{T1}, \cite{T2} and \cite{T3} for more details on
this topic. In our particular case $X$ is a distance-regular graph
with $\partial (x,y) = |x \triangle y|$ (where $\triangle$ stands
for symmetric difference) and $X_r$ is a distance-regular graph
with $\partial_r (x,y) = |x \setminus y|$ (Hamming and Johnson
graphs). When $Y = X$ and $x_0 = \emptyset$, the corresponding
Terwilliger algebra is $\mbox{End}_G(L^2(X))$, the algebra we are
interested in. The algebra $\mbox{End}_G(L^2(X))$ connects the
Bose-Mesner algebras of Johnson graphs for different values of
$r$.

We would like to point out that, after this paper was submitted
for pu\-blication, the referee communicated to us the existence of
Go's article \cite{Go}, which is related to the present paper. The
central topic of Go's paper is the Terwilliger algebra of the
hypercube, $X$ above. That is, in \cite{Go} the algebra
$\mbox{End}_G(L^2(X))$ is studied from a different point of view.
She regards this algebra as a homomorphic image of the universal
enveloping algebra of $\mathfrak{sl}(2,\mathbb{C})$; and then she
works with operators defined in terms of two natural generators
$A, A^{\star}$ with kernels $$a(x,y) = \left\{ \begin{array}{ll} 1
& \mbox{if} \ \ |x \triangle y| = 1 \\ 0 & \mbox{otherwise}
\end{array} \right., \qquad a^{\star}(x,y) = \left\{
\begin{array}{ll} n-2|x| & \mbox{if} \ \ x = y \\ 0 &
\mbox{otherwise.} \end{array} \right.$$ She studies the
irreducible submodules and obtains expressions for the central
primi\-tive idempotents of the algebra. We consider this is an
interesting approach, but it is not easy (unless one introduces
further arguments) to obtain our results using the information
contained in Go's paper.

The organization of our paper is as follows. In section
$\ref{section2}$ we define the notions of finite symmetric space
and spherical function, and then we recall some basic results that
are used all throughout the paper. In section $\ref{section3}$ we
give the decomposition of $\rho$ into irreducible components and
we analyze the kernels $\la$ via the Radon transforms. In section
$\ref{section4}$ we introduce a discrete Laplacian operator and
then we show how the kernels $\la$ can be viewed as solutions of a
hypergeometric equation. We study that equation in detail.
Finally, in section \ref{section5}, we deal with the mentioned
analysis of the Fourier transform $\mathcal{F}_X$.

\section{Finite symmetric spaces}
\label{section2}

We begin with a summary of some basic results about finite
symme\-tric spaces and spherical functions that will be used in
the sequel. For further information on these topics see \cite{Te}
and the references cited there. Let $G$ be a finite group acting
on a finite set $X$, this action leads us to the associated
unitary representation $$\rho: G \longrightarrow U(L^2(X))$$ given
by $\big( \rho(g) \psi \big) (x) = \psi (g^{-1} x)$. Assume the
action is transitive, then $X$ is said to be a \emph{finite
symmetric space} with respect to $G$ if the algebra $\mbox{End}_G
(L^2(X))$ of the endomorphisms on $L^2(X)$ which commute with the
action of $\rho$ is abelian.

\begin{Rem}
\emph{We recall that $\mbox{End}_G (L^2(X))$ is an abelian algebra
if and only if the representation $\rho$ is multiplicity-free. So
we can invoke this classical result of representation theory to
give another characterization of finite symmetric spaces.}
\end{Rem}
Now assume we are given a couple of finite symmetric spaces $X_1$
and $X_2$ with respect to $G$. Let us denote by $\rho_1$ and
$\rho_2$ the respective associated representations. We assign to
each operator $T \in \mbox{Hom} (L^2(X_1),L^2(X_2))$ the matrix
$\xi$ of $T$ with res\-pect to the natural bases of $L^2(X_1)$ and
$L^2(X_2)$. This mapping is clearly a linear isomorphism from
$\mbox{Hom} (L^2(X_1),L^2(X_2))$ onto $L^2(X_2 \times X_1)$, we
denote it by $\Psi$. Thus $T$ and $\xi$ are related by the
expression $$(T \psi)(x_2) = \sum_{x_1 \in X_1} \xi(x_2,x_1)
\psi(x_1).$$

If we compare the operators $T \circ \rho_1(g)$ and $\rho_2(g)
\circ T$ written in this way, it is obvious that $T \in
\mbox{Hom}_G (L^2(X_1),L^2(X_2))$ if and only if the relation
$\xi(gx_2,gx_1) = \xi(x_2,x_1)$ holds for all $(x_2,x_1) \in X_2
\times X_1$ and all $g \in G$. Here $\mbox{Hom}_G
(L^2(X_1),L^2(X_2))$ denotes the algebra of intertwining operators
for $\rho_1$ and $\rho_2$. That is, $T$ is an intertwining
operator for $\rho_1$ and $\rho_2$ if and only if the associated
matrix is constant at the orbits of the action
$$\begin{array}{rcl} G \times X_2 \times X_1 & \longrightarrow &
X_2 \times X_1 \\ (g, (x_2,x_1)) & \longmapsto & (g x_2, g x_1).
\end{array}$$

An action of $G$ on a finite set $X$ is called \emph{symmetric} if
for all $x,x' \in X$ there exists $g \in G$ such that $g x = x'$
and $g x' = x$. A finite set $X$ endowed with a symmetric action
of $G$ is automatically a finite symmetric space with respect to
$G$. To justify this we observe that if the action of $G$ on $X$
is symmetric then $\Psi \big( \mbox{End}_G(L^2(X) \big)$ is a
subalgebra of $L^2(X \times X)$ made up of symmetric matrices,
hence abelian. Now, taking into account that $\Psi$ is an algebra
isomorphism when $X_1 = X_2$, the result follows. Let us consider
the set $\G_X = \{\pi \in \G: \mbox{Mult}_{\pi} (\rho) \neq 0\}$,
where $\G$ stands for the dual object, the set of irreducible
unitary representations of $G$. Note that if $X$ is symmetric with
respect to $G$, then every $\pi \in \G_X$ satisfies
$\mbox{Mult}_{\pi} (\rho) = 1$, since $\rho$ is multiplicity-free.
Then we use the set $\G_X$ to decompose the space $L^2(X)$ into
irreducible components $$L^2(X) = \bigoplus_{\pi \in \G_X}
L^2(X)_{\pi}.$$ We write $P_{\pi}$ for the orthogonal projection
onto $L^2(X)_{\pi}$, and the matrix of $P_{\pi}$ will be denoted
by $p_{\pi}$. The \emph{spherical functions} on $X$ are defined by
$$\xi_{X,\pi} = \frac{|X|}{d(\pi)} p_{\pi} \in \Psi\big(
\mbox{End}_G (L^2(X)) \big)$$ where $\pi \in \G_X$ and $d(\pi)$
denotes the degree of $\pi$. We also write $\mathcal{S}_{X,\pi}$
for the \emph{associated operator} in $\textnormal{End}_G(L^2(X))$
with matrix $\xi_{X,\pi}$. The proof of the following theorem can
be found in \cite{Te}.

\begin{The} \label{theorem1}
Let $X$ be a finite symmetric space with respect to the finite
group $G$ and let $\xi \in \Psi\big( \mbox{End}_G (L^2(X)) \big)$,
then the following are equi\-valent:
\begin{itemize}
\item[$(a)$] There exists $\pi \in \G_X$ such that $\xi =
\xi_{X,\pi}$.
\item[$(b)$] $\ \xi(x_0,x_0) = 1$ for all $x_0 \in X$ and for every
$x_1,x_2 \in X$ $$\frac{1}{|G_{x_0}|} \sum_{g \in G_{x_0}} \xi(g
x_1,x_2) = \xi(x_1,x_0) \xi(x_0,x_2)$$ where $G_{x_0}$ denotes the
isotropy subgroup of $x_0$.
\end{itemize}
\end{The}

\section{The algebra $\textnormal{E\lowercase{nd}}_{S(\Omega)} \big(
L^2(\mathcal{P}(\Omega)) \big)$} \label{section3}

\setcounter{The}{0}

As we pointed out in the introduction, the symmetric group
$S(\Omega)$ acts naturally on the power set $\mathcal{P} (\Omega)$
providing the associated unitary representation $\rho$. We recall
that $G$ stands for $S(\Omega)$ and $X$ for $\mathcal{P}
(\Omega)$. $X$ is not a symmetric space with respect to $G$, in
fact the mentioned action is not even transitive. Nevertheless the
orbits of such action are given by the family of sets $X_r = \{x
\in X: |x|=r\}$, $0 \le r \le n$. This action is symmetric on each
orbit, so we know that the sets $X_r$ are symmetric spaces with
respect to $G$ for $0 \le r \le n$. If we denote by $\rho_r: G
\longrightarrow U ( L^2(X_r))$ the associa\-ted representations
and we identify the space $L^2(X_r)$ with the subspace of $L^2(X)$
of functions supported on $X_r$, then it is very easy to check
that $$\rho = \bigoplus_{0 \le r \le n} \rho_r.$$

We recall that the matrix of an operator $T \in \mbox{Hom}_{G}
(L^2(X_{r_1}),L^2(X_{r_2}))$ is cons\-tant at the orbits $\theta_k
= \big\{ (x_2,x_1) \in X_{r_2} \times X_{r_1}: |x_2 \setminus x_1|
= k \big\}$ of the natural action of $G$ on $X_{r_2} \times
X_{r_1}$. Here $0 \vee (r_2 - r_1) \le k \le (n - r_1) \wedge r_2$
and $k \in \N$. So we can write these operators in the form $$(T
\psi) (x_2) = \sum_{x_1 \in X_{r_1}} \xi(|x_2 \setminus x_1|)
\psi(x_1)$$ where $\xi$ depends on the variable $k$. The function
$\xi$ is called the \emph{kernel} of $T$. We also know that the
dimension of $\mbox{Hom}_{G} (L^2(X_{r_1}),L^2(X_{r_2}))$
coincides with the number of orbits $\theta_k$. That is $$\dim
\big( \mbox{Hom}_{G} (L^2(X_{r_1}),L^2(X_{r_2})) \big) = r_1
\wedge (n - r_1) \wedge r_2 \wedge (n - r_2) + 1.$$ Now, taking
into account that $X_r$ is a finite symmetric space with respect
to $G$, we deduce that the spaces $L^2(X_r)$ are
multiplicity-free. Therefore, according to Schur's lemma, the
dimension of $\mbox{Hom}_{G} (L^2(X_{r_1}),L^2(X_{r_2}))$ gives
the number of irreducible components that $L^2(X_{r_1})$ and
$L^2(X_{r_2})$ have in common. In particular, if $[n/2]$ denotes
the integer part of $n/2$
\begin{enumerate}
\item For $0 \le r \le [n/2]$, the space $L^2(X_r)$ has $r + 1$
irreducible components.
\item For $0 \le r < [n/2]$, the spaces $L^2(X_r)$ and
$L^2(X_{r+1})$ have $r + 1$ irredu\-cible components in common.
\end{enumerate}
Hence, by a simple induction argument, there exist a family of
inequivalent irreducible representations $\pi_s: G \longrightarrow
U(\mathcal{V}_s)$ where $0 \le s \le [n/2]$ and such that $$\rho_r
\simeq \bigoplus_{0 \le s \le r} \pi_s$$ for $0 \le r \le [n/2]$.
On the other hand the representations $\rho_r$ and $\rho_{n-r}$
are equivalent. Namely, the operator $\mathcal{C}^r: L^2(X_r)
\longrightarrow L^2(X_{n-r})$ defined by $(\mathcal{C}^r \psi)
(x^c) = \psi(x)$ is an intertwining unitary operator. Thus, for $0
\le r \le n$, we have $$L^2(X_r) \simeq \bigoplus_{s = 0}^{r
\wedge (n-r)} \mathcal{V}_s.$$ We shall denote by $L^2(X_r)_s$ the
$G$-invariant subspace of $L^2(X_r)$ equivalent to
$\mathcal{V}_s$. Finally we note that
\begin{equation} \label{equation1}
\dim (\mathcal{V}_s) = \dim (L^2(X_s)) - \dim (L^2(X_{s-1})) =
{{n}\choose{s}} - {{n}\choose{s-1}}
\end{equation}

\begin{Rem}
\emph{The study of representations of the symmetric group provides
techniques, such as the Young's rule, that can be used to
determine how are the representations $\pi_s$ for $0 \le s \le
[n/2]$. The result is that $\pi_s$ coincides with the irreducible
representation $\pi_{(n-s,s)}$ associated to the arithmetic
partition $(n-s,s)$ of $n$. See \cite{FH} for the details.}
\end{Rem}
Once we know the irreducible components of $\rho$, we have the
following decomposition for the algebra $\mbox{End}_{G} (L^2(X))$
\begin{eqnarray*}
\mbox{End}_{G} (L^2(X)) & = & \bigoplus_{0 \le r_1,r_2 \le n}
\mbox{Hom}_{G} (L^2(X_{r_1}),L^2(X_{r_2})) \\ & = & \bigoplus_{0
\le r_1,r_2 \le n} \bigoplus_{s=0}^{N(r_1,r_2)} \mbox{Hom}_{G}
\big( L^2(X_{r_1})_s,L^2(X_{r_2})_s \big)
\end{eqnarray*}
with $N(r_1,r_2) = r_1 \wedge (n-r_1) \wedge r_2 \wedge (n-r_2)$.
Now, by Schur's lemma, all the spaces $\mbox{Hom}_{G} \big(
L^2(X_{r_1})_s,L^2(X_{r_2})_s \big)$ are $1$-dimensional.

\begin{Def} \emph{Let $0 \le r_1,r_2 \le n$ and $0 \le s \le
N(r_1,r_2)$. We define the \emph{operator $\La$} as a non zero
element of $\mbox{Hom}_{G} \big( L^2(X_{r_1})_s,L^2(X_{r_2})_s
\big)$.}
\end{Def}

\begin{Rem}
\emph{The definition of the operators $\La$ is ambiguous, we will
normalize these operators after lemma \ref{lemma1}. Note that,
defining $\La$ by 0 on the subspace orthogonal to
$L^2(X_{r_1})_s$, these operators are elements of the algebra
$\mbox{End}_{G} (L^2(X))$. We also note that $\Lambda_s^{r_1,r_2}
\circ \Lambda_{s'}^{r_3,r_4} = 0$ unless $s = s'$ and $r_1 =
r_4$.}
\end{Rem}

\begin{Pro}
The family of operators $\La$ is an orthogonal basis of the space
$\textnormal{End}_{G} (L^2(X))$ with respect to the
Hilbert-Schmidt inner product.
\end{Pro}

\begin{proof}
This family is obviously a basis of $\mbox{End}_G (L^2(X))$. To
see the orthogonality we observe that, since the adjoint of an
intertwining operator is also an intertwining operator, there
exist a non zero constant $c_s(r_1,r_2)$ such that $(\La)^{\star}
= c_s(r_1,r_2) \Lambda_s^{r_2,r_1}$. Hence we can write
$$\mbox{tr} \big( \La \circ (\Lambda_{s'}^{r_3,r_4})^{\star} \big)
= c_{s'}(r_3,r_4)  \,\ \mbox{tr} \big( \La \circ
\Lambda_{s'}^{r_4,r_3} \big)$$ which is 0 unless $s = s'$,
$r_1=r_3$ and $r_2 = r_4$. This completes the proof.
\end{proof}

We now define the \emph{Radon transforms} $\Ra_{\subset}^{r_1,r_2}
\in \mbox{Hom}_G (L^2(X_{r_1}),L^2(X_{r_2}))$ for $r_1 \le r_2$
and their adjoints $(\Ra_{\subset}^{r_1,r_2})^{\star} \in
\mbox{Hom}_G (L^2(X_{r_2}),L^2(X_{r_1}))$ by
$$(\Ra_{\subset}^{r_1,r_2} \psi) (x_2) = \sum_{x_1 \subset x_2}
\psi(x_1) \qquad (\Ra_{\supset}^{r_2,r_1} \psi) (x_1) = \sum_{x_2
\supset x_1} \psi(x_2).$$ We know that $\La$ is an intertwining
operator, thus we can write its kernel $\lambda_s^{r_1,r_2}$ as a
function of the integer variable $k$, where we recall that $$0
\vee (r_2-r_1) \le k \le (n-r_1) \wedge r_2.$$ The Radon
transforms will be very useful in the study of the kernels $\la$.
Given $\alpha \in \R$ and $k \in \N$ we recall the classical
notation $[\alpha]_k = \alpha (\alpha - 1) \cdots (\alpha - k +
1)$ $([\alpha]_0 = 1)$ and $(\alpha)_k = \alpha (\alpha + 1)
\cdots (\alpha + k - 1)$ $((\alpha)_0 = 1)$.

\begin{Lem}
Given $0 \le r \le n$, $0 \le s \le r \wedge (n-r)$ and $0 \le k
\le s$, we have $$\lambda_s^{r,s} (k) = (-1)^k
\frac{(r-s+1)_k}{[n-r]_k} \lambda_s^{r,s} (0).$$
\end{Lem}

\begin{proof}
The relation is obvious for $s = 0$. Otherwise we observe that
$\mbox{Im} (\Lambda_s^{r,s}) = L^2(X_s)_s \subset \mbox{Ker}
(\Ra_{\supset}^{s,s-1})$. In particular $\Ra_{\supset}^{s,s-1}
\circ \Lambda_s^{r,s} = 0$. Then, if we define $\delta_x \in
L^2(X_r)$ to be $1$ at $x \in X_r$ and $0$ otherwise we take $y
\in X_{s-1}$ to get
\begin{eqnarray*} \big( (\Ra_{\supset}^{s,s-1} \circ
\Lambda_s^{r,s}) \delta_x \big) (y) & = & \sum_{z \supset y}
(\Lambda_s^{r,s} \delta_x) (z) = \sum_{z \supset y}
\lambda_s^{r,s} (|z \setminus x|) \\ & = & |x \setminus y|
\lambda_s^{r,s} (|y \setminus x|) + |(\Omega \setminus x)
\setminus y| \lambda_s^{r,s} (|y \setminus x| + 1).
\end{eqnarray*}
Therefore, we have the following recurrence equation for $1 \le k
\le s$ $$(r-s+k) \lambda_s^{r,s} (k-1) + (n-r-k+1) \lambda_s^{r,s}
(k) = 0.$$ Solving the recurrence equation we get the desired
relation.
\end{proof}

Now we want to show that there exists a constant $C_s(r_1,r_2)$
such that the fo\-llowing relation holds
\begin{equation} \label{equation2}
\la (k) = C_s(r_1,r_2) \sum_{j=0}^s {{k}\choose{j}}
{{r_2-k}\choose{s-j}} \lambda_s^{r_1,s} (j)
\end{equation}
in the usual rank $0 \le r_1,r_2 \le n$, $0 \le s \le r_1 \wedge
(n-r_1) \wedge r_2 \wedge (n-r_2)$ and $0 \vee (r_2 - r_1) \le k
\le (n-r_1) \wedge r_2$. To prove this we observe that, since the
Radon transforms are intertwining operators, there exists a
constant $c_s^{r_1,r_2}$ such that $\Ra_{\subset}^{s,r_2} \circ
\Lambda_s^{r_1,s} = c_s^{r_1,r_2} \La$. Now, if $x_1 \in X_{r_1}$
and $x_2 \in X_{r_2}$, then
\begin{eqnarray*}
c_s^{r_1,r_2} \la (|x_2 \setminus x_1|) & = & \big(
(\Ra_{\subset}^{s,r_2} \circ \Lambda_s^{r_1,s}) \delta_{x_1} \big)
(x_2) = \sum_{z \subset x_2} \lambda_s^{r_1,s} (|z \setminus x_1|)
\\ & = & \sum_{j=0}^s \big| \big\{ z \in X_s: \ \ z \subset x_2, \
\ |z \setminus x_1| = j \big\} \big| \ \ \lambda_s^{r_1,s} (j)
\\ & = & \sum_{j=0}^s {{|x_2 \setminus x_1|}\choose{j}} {{r_2 -
|x_2 \setminus x_1|}\choose{s-j}} \lambda_s^{r_1,s}(j).
\end{eqnarray*}
Taking $k = |x_2 \setminus x_1|$ the relation in (\ref{equation2})
arises.

\begin{Rem}
\emph{Relation (\ref{equation2}) gives that $\la$ is a polynomial
of degree $\le s$. But its highest coefficient is a non zero
factor of $$\sum_{j=0}^s \frac{(-1)^{s-j}}{j!(s-j)!}
\lambda_s^{r_1,s} (j) = \lambda_s^{r_1,s} (0) \sum_{j=0}^s
\frac{(-1)^s}{j!(s-j)!} \frac{(r_1-s+1)_j}{[n-r_1]_j} \neq 0$$ and
so we deduce that $\partial \la = s$ and $\la (0) \neq 0$. On the
other hand the polynomial $\la$ is defined at $r_1 \wedge (n-r_1)
\wedge r_2 \wedge (n-r_2) +1$ points and this number is always
greater than $s$. Thus we know that there exists a unique
polynomial $p \in \C [t]$ of degree $s$ which extends $\la$. In
what follows we shall denote that polynomial with the same
expression $\la$.}
\end{Rem}
Using the remark above we can evaluate equation (\ref{equation2})
at $t = 0$ to obtain the following result.
\begin{Lem} \label{lemma1}
Given $r_1,r_2$ and $s$ in the usual rank of parameters and $t \in
\C$, we have $$\frac{\la (t)}{\la (0)} = {{r_2}\choose{s}}^{-1}
\sum_{j=0}^s {{t}\choose{j}} {{r_2-t}\choose{s-j}}
\frac{\lambda_s^{r_1,s} (j)}{\lambda_s^{r_1,s} (0)}.$$
\end{Lem}
We are now in a position to normalize the operators $\La$. Just
take $\la (0)$ to be 1 for all possible values of $r_1,r_2$ and
$s$.
\begin{Rem}
\emph{We shall write $\Lambda_s^r$ and $\lambda_s^r $ for
$\Lambda_s^{r,r}$ and $\lambda_s^{r,r}$. With this
norma\-lization, the mappings $$(x,x') \in X_r \times X_r
\longmapsto \lambda_s^r (|x \setminus x'|) \in \R$$ with $0 \le r
\le n$ and $0 \le s \le r \wedge (n-r)$, are the spherical
functions on $X_r$.}
\end{Rem}
We now want to find the value of the kernels $\la$ at $t = r_2$.
For that we just need to evaluate the expression given in lemma
\ref{lemma1} at $t = r_2$
\begin{equation} \label{equation3}
\la (r_2) = \lambda_s^{r_1,s} (s) = (-1)^s
\frac{[r_1]_s}{[n-r_1]_s}.
\end{equation}
In the following result we use identity (\ref{equation3}) to
relate the operator $\La$ with some other operators of the same
basis.

\begin{Pro} \label{proposition1}
Given $r_1,r_2$ and $s$ in the usual rank of parameters, we have
\begin{eqnarray*}
\La \circ \mathcal{C}^{n-r_1} & = & (-1)^s
\frac{[r_1]_s}{[n-r_1]_s} \Lambda_s^{n-r_1,r_2} \\
\mathcal{C}^{r_2} \circ \La & = & (-1)^s \frac{[n-r_2]_s}{[r_2]_s}
\Lambda_s^{r_1,n-r_2} \\ (\La)^{\star} & = &
\frac{[r_1]_s[n-r_2]_s}{[n-r_1]_s[r_2]_s} \Lambda_s^{r_2,r_1}.
\end{eqnarray*}
\end{Pro}

\begin{proof}
The operators $\mathcal{C}^r$ belong to the algebra
$\mbox{End}_{G} (L^2(X))$, so there exist cons\-tants
$c_s^1(r_1,r_2), c_s^2(r_1,r_2)$ and $c_s^3(r_1,r_2)$ such that
\begin{eqnarray*}
\La \circ \mathcal{C}^{n-r_1} & = & c_s^1(r_1,r_2)
\Lambda_s^{n-r_1,r_2} \\ \mathcal{C}^{r_2} \circ \La & = &
c_s^2(r_1,r_2) \Lambda_s^{r_1,n-r_2} \\ \La & = & c_s^3(r_1,r_2)
(\Lambda_s^{r_2,r_1})^{\star}.
\end{eqnarray*}
If we express these equalities in terms of the kernels we get
\begin{eqnarray*}
\la (t) & = & c_s^1(r_1,r_2) \lambda_s^{n-r_1,r_2} (r_2 - t) \\
\la (t) & = & c_s^2(r_1,r_2) \lambda_s^{r_1,n-r_2} (n-r_1-t) \\
\la (t) & = & c_s^3(r_1,r_2) \lambda_s^{r_2,r_1} (r_1-r_2+t).
\end{eqnarray*}
Taking $t = r_2$ in the first equation and using (\ref{equation3})
we obtain the value of $c_s^1(r_1,r_2)$. The same occurs in the
third equation. Finally, if we evaluate the third equation at $t =
r_2-r_1$ we get $$\la (r_2-r_1) =
\frac{[r_1]_s[n-r_2]_s}{[n-r_1]_s[r_2]_s}.$$ We conclude by
evaluating the second equation at $t = r_2-r_1$.
\end{proof}

The definition of $\La$ gives $\Lambda_s^{r_2,r_3} \circ \La =
C_s(r_1,r_2,r_3) \Lambda_s^{r_1,r_3}$ in the usual rank of
parameters. In the following theorem we investigate the value of
the cocycle $C_s(r_1,r_2,r_3)$, obtaining in such a way a
\emph{multiplication table} for the operators $\La$.

\begin{The} \label{theorem2}
The operators $\La$ satisfy, in the usual rank of para\-meters,
the following relation $$\Lambda_s^{r_2,r_3} \circ \La =
\frac{{{n}\choose{r_2}}}{{{n}\choose{s}} - {{n}\choose{s-1}}}
\Lambda_s^{r_1,r_3}.$$
\end{The}

\begin{proof}
We recall that $\dim \big( L^2(X_r)_s \big) = \dim (\mathcal{V}_s)
= {{n}\choose{s}} - {{n}\choose{s-1}}$, see equation
(\ref{equation1}). We also note that $$\mbox{tr} (\Lambda_s^r) =
\sum_{x \in X_r} \lambda_s^r (|x \setminus x|) =
{{n}\choose{r}}.$$ Theses calculations give us an explicit form of
the orthogonal projection of $L^2(X)$ onto $L^2(X_r)_s$ $$P_s^r =
\frac{{{n}\choose{s}} - {{n}\choose{s-1}}}{{{n}\choose{r}}}
\Lambda_s^r.$$ Hence we can write
\begin{equation} \label{equation4}
\La \circ \Lambda_s^{r_1} =
\frac{{{n}\choose{r_1}}}{{{n}\choose{s}} - {{n}\choose{s-1}}} \La
\ \ ; \ \ \Lambda_s^{r_2} \circ \La =
\frac{{{n}\choose{r_2}}}{{{n}\choose{s}} - {{n}\choose{s-1}}} \La
\end{equation}
We then have, in terms of the kernels, the relations
\begin{equation} \label{equation5}
\begin{array}{c}
\displaystyle \sum_{y \in X_{r_1}} \la (|x_2 \setminus y|)
\lambda_s^{r_1} (|y \setminus x_1|) =
\frac{{{n}\choose{r_1}}}{{{n}\choose{s}} - {{n}\choose{s-1}}} \la
(|x_2 \setminus x_1|) \\ \displaystyle \sum_{y \in X_{r_2}}
\lambda_s^{r_2} (|x_2 \setminus y|) \la (|y \setminus x_1|) =
\frac{{{n}\choose{r_2}}}{{{n}\choose{s}} - {{n}\choose{s-1}}} \la
(|x_2 \setminus x_1|)
\end{array}
\end{equation}
We recall that the mappings $(x,x') \in X_r \times X_r \longmapsto
\lambda_s^r(|x \setminus x'|) \in \R$ are the spherical functions
associated to the symmetric spaces $X_r$. Thus we invoke theorem
\ref{theorem1} and (\ref{equation5}) to get the following
relations where $x_0 \in X_{r_1}$ and $x_0' \in X_{r_2}$
\begin{equation} \label{equation6}
\begin{array}{c} \displaystyle
\frac{1}{|G_{x_0}|} \sum_{g \in G_{x_0}} \la (|x_2 \setminus g
x_1|) = \la (|x_2 \setminus x_0|) \lambda_s^{r_1} (|x_0 \setminus
x_1|) \\ \displaystyle \frac{1}{|G_{x'_0}|} \sum_{g \in G_{x'_0}}
\la (|g x_2 \setminus x_1|) = \lambda_s^{r_2} (|x_2 \setminus
x'_0|) \la (|x'_0 \setminus x_1|)
\end{array}
\end{equation}
We now fix $x_j \in X_{r_j}$ with $j=1,2,3$. Let us assume for the
moment that $r_1 \ge r_2 \ge r_3$ and $x_3 \subset x_2 \subset
x_1$, then
\begin{eqnarray*}
C_s(r_1,r_2,r_3) & = & \frac{C_s(r_1,r_2,r_3)}{|G_{x_2}|} \sum_{g
\in G_{x_2}} \lambda_s^{r_1,r_3}(|g x_3 \setminus x_1|) \\ & = &
\frac{1}{|G_{x_2}|} \sum_{g \in G_{x_2}} \sum_{y \in X_{r_2}}
\lambda_s^{r_2,r_3} (|g x_3 \setminus y|) \la (|y \setminus x_1|)
\\ & = & \sum_{y \in X_{r_2}} \Big( \frac{1}{|G_{x_2}|}
\sum_{g \in G_{x_2}} \lambda_s^{r_2,r_3} (|x_3 \setminus g y|)
\Big) \la (|y \setminus x_1|)
\end{eqnarray*}
and using (\ref{equation5}) and (\ref{equation6}) we obtain
\begin{eqnarray*}
C_s(r_1,r_2,r_3) & = & \lambda_s^{r_2,r_3} (|x_3 \setminus x_2|)
\sum_{y \in X_{r_2}} \lambda_s^{r_2} (|x_2 \setminus y|) \la (|y
\setminus x_1|) \\ & = & {{n}\choose{r_2}}
\frac{\lambda_s^{r_2,r_3} (|x_3 \setminus x_2|) \la (|x_2
\setminus x_1|)}{{{n}\choose{s}} - {{n}\choose{s-1}}} =
\frac{{{n}\choose{r_2}}}{{{n}\choose{s}} - {{n}\choose{s-1}}}.
\end{eqnarray*}

Therefore we have proved the desired relation when $r_1 \ge r_2
\ge r_3$ and by (\ref{equation4}) also for $r_1 = r_2$ and $r_2 =
r_3$. But we want to know the value of $C_s(r_1,r_2,r_3)$ for all
possible values of $r_1,r_2$ and $r_3$ in
\begin{equation} \label{equation7}
\Lambda_s^{r_2,r_3} \circ \La = C_s(r_1,r_2,r_3)
\Lambda_s^{r_1,r_3}
\end{equation}
Taking adjoints in (\ref{equation7}) and applying proposition
\ref{proposition1} we get
\begin{equation} \label{equation8}
C_s(r_1,r_2,r_3) = C_s(r_3,r_2,r_1).
\end{equation}
Now we compose with the suitable $\mathcal{C}^r$ operators
\begin{eqnarray*}
\Lambda_s^{r_2,r_3} \circ \La \circ \mathcal{C}^{n-r_1} & = &
C_s(r_1,r_2,r_3) \ \ \Lambda_s^{r_1,r_3} \circ \mathcal{C}^{n-r_1}
\\ \mathcal{C}^{r_3} \circ \Lambda_s^{r_2,r_3} \circ \La & = &
C_s(r_1,r_2,r_3) \ \ \mathcal{C}^{r_3} \circ \Lambda_s^{r_1,r_3}
\\ \Lambda_s^{r_2,r_3} \circ \mathcal{C}^{n-r_2} \circ
\mathcal{C}^{r_2} \circ \La & = & C_s(r_1,r_2,r_3) \ \
\Lambda_s^{r_1,r_3}
\end{eqnarray*}
and then proposition \ref{proposition1} gives rise to
\begin{eqnarray}
\nonumber C_s(r_1,r_2,r_3) & = & C_s(n-r_1,r_2,r_3) \\ & = &
C_s(r_1,n-r_2,r_3) \label{equation9} \\ \nonumber & = &
C_s(r_1,r_2,n-r_3).
\end{eqnarray}
Equation (\ref{equation8}) reduces the possibilities to the case
$r_1 \ge r_3$. But we know the value of $C_s(r_1,r_2,r_3)$ for
$r_1 \ge r_2 \ge r_3$, thus we only consider the cases $r_1 \ge
r_3 \ge r_2$ and $r_2 \ge r_1 \ge r_3$. In the first case we know
that $C_s(n,r_1,r_2) = C_s(n,r_1,r_3)$, so right multiplication by
$\Lambda_s^{n,r_1}$ in (\ref{equation7}) implies that
$$\Lambda_s^{r_2,r_3} \circ \Lambda_s^{n,r_2} = C_s(r_1,r_2,r_3)
\Lambda_s^{n,r_3}.$$ Equalities (\ref{equation9}) give
$\Lambda_s^{n-r_2,n-r_3} \circ \Lambda_s^{n,n-r_2} =
C_s(r_1,r_2,r_3) \Lambda_s^{n,n-r_3}$. Therefore
$$C_s(r_1,r_2,r_3) = C_s(n,n-r_2,n-r_3) =
\frac{{{n}\choose{n-r_2}}}{{{n}\choose{s}} - {{n}\choose{s-1}}} =
\frac{{{n}\choose{r_2}}}{{{n}\choose{s}} - {{n}\choose{s-1}}}.$$
In the second case $C_s(r_2,r_3,0) = C_s(r_1,r_3,0)$, hence left
multiplication by $\Lambda_s^{r_3,0}$ in (\ref{equation7}) gives
$\Lambda_s^{r_2,0} \circ \La = C_s(r_1,r_2,r_3)
\Lambda_s^{r_1,0}$. Now by (\ref{equation9}) we get
$$\Lambda_s^{n-r_2,0} \circ \Lambda_s^{n-r_1,n-r_2} =
C_s(r_1,r_2,r_3) \Lambda_s^{n-r_1,0}.$$ So we have
$$C_s(r_1,r_2,r_3) = C_s(n-r_1,n-r_2,0) =
\frac{{{n}\choose{n-r_2}}}{{{n}\choose{s}} - {{n}\choose{s-1}}} =
\frac{{{n}\choose{r_2}}}{{{n}\choose{s}} - {{n}\choose{s-1}}}.$$
This concludes the proof.
\end{proof}

Making use of theorem \ref{theorem2} we can compute the
Hilbert-Schmidt norm of the operators $\La$.

\begin{Cor} \label{corollary1}
Given $r_1,r_2$ and $s$ in the usual rank of parameters, we have
$$\textnormal{tr} \big( \La \circ (\La)^{\star} \big) =
\frac{[r_1]_s[n-r_2]_s}{[n-r_1]_s[r_2]_s}
\frac{{{n}\choose{r_1}}{{n}\choose{r_2}}}{{{n}\choose{s}} -
{{n}\choose{s-1}}}.$$
\end{Cor}

\begin{proof}
We just apply proposition \ref{proposition1} and theorem
\ref{theorem2} consecutively
\begin{eqnarray*}
\mbox{tr} \big( \La \circ (\La)^{\star} \big) & = &
\frac{[r_1]_s[n-r_2]_s}{[n-r_1]_s[r_2]_s} \mbox{tr} \big( \La
\circ \Lambda_s^{r_2,r_1} \big) \\ & = &
\frac{[r_1]_s[n-r_2]_s}{[n-r_1]_s[r_2]_s}
\frac{{{n}\choose{r_1}}}{{{n}\choose{s}} - {{n}\choose{s-1}}} \ \
\mbox{tr} (\Lambda_s^{r_2}) \\ & = &
\frac{[r_1]_s[n-r_2]_s}{[n-r_1]_s[r_2]_s}
\frac{{{n}\choose{r_1}}{{n}\choose{r_2}}}{{{n}\choose{s}} -
{{n}\choose{s-1}}}.
\end{eqnarray*}
This concludes the proof.
\end{proof}

\begin{Rem}
\emph{The normalization we have chosen for the operators $\La$ is
very natural since we have obtained from it the spherical
functions on the symmetric spaces $X_r$. Anyway there are other
ways to normalize the ope\-rators $\La$ which are also very
significant, for instance the normalization
\begin{equation} \label{equation10}
\widetilde{\Lambda}_s^{r_1,r_2} =
\frac{[n-r_1]_s[r_2]_s}{\sqrt{{{n-2s}\choose{r_1-s}}
{{n-2s}\choose{r_2-s}}} s![n-s+1]_s} \La
\end{equation}
provides the simpler relations $$\begin{array}{rclcrcl}
\widetilde{\Lambda}_s^{r_1,r_2} \circ \mathcal{C}^{n-r_1} & = &
(-1)^s \widetilde{\Lambda}_s^{n-r_1,r_2} & \qquad &
(\widetilde{\Lambda}_s^{r_1,r_2})^{\star}  & = &
\widetilde{\Lambda}_s^{r_2,r_1}
\\ \mathcal{C}^{r_2} \circ
\widetilde{\Lambda}_s^{r_1,r_2} & = & (-1)^s
\widetilde{\Lambda}_s^{r_1,n-r_2} & \qquad &
\widetilde{\Lambda}_s^{r_2,r_3} \circ
\widetilde{\Lambda}_s^{r_1,r_2} & = &
\widetilde{\Lambda}_s^{r_1,r_2}.
\end{array}$$ Moreover, the operators
$\widetilde{\Lambda}_s^{r_1,r_2}$ are unitary isomorphisms from
$L^2(X_{r_1})_s$ onto the space $L^2(X_{r_2})_s$. We shall denote
the normalization constant by $\alpha_s^{r_1,r_2}$. Also the
normalization $$\overline{\Lambda}_s^{r_1,r_2} =
\frac{1}{\sqrt{{{n}\choose{s}} - {{n}\choose{s-1}}}}
\widetilde{\Lambda}_s^{r_1,r_2}$$ gives an orthonormal basis with
respect to the Hilbert-Schmidt inner pro\-duct. With this
normalization we have the relations $$\begin{array}{rclcrcl}
\overline{\Lambda}_s^{r_1,r_2} \circ \mathcal{C}^{n-r_1} & = &
(-1)^s \overline{\Lambda}_s^{n-r_1,r_2} & \qquad &
(\overline{\Lambda}_s^{r_1,r_2})^{\star}  & = &
\overline{\Lambda}_s^{r_2,r_1}
\\ \mathcal{C}^{r_2} \circ
\overline{\Lambda}_s^{r_1,r_2} & = & (-1)^s
\overline{\Lambda}_s^{r_1,n-r_2} & \qquad &
\overline{\Lambda}_s^{r_2,r_3} \circ
\overline{\Lambda}_s^{r_1,r_2} & = & \beta_s
\overline{\Lambda}_s^{r_1,r_2}.
\end{array}$$ where $\beta_s$ denotes the normalization constant
$\Big( \sqrt{{{n}\choose{s}} - {{n}\choose{s-1}}} \Big)^{-1}$.}
\end{Rem}

\section{The use of a Laplacian operator}
\label{section4}

\setcounter{The}{0}

We now investigate some useful expressions for the kernels $\la$
by showing the associated operators $\La$ as eigenvectors of a
Laplacian operator. For that purpose it will be necessary to
understand the cardinal of $\Omega$ as another parameter of the
problem. That is, the number $n$ is no longer a fixed value and so
we will often explicit the dependence on $n$ for clarity. We start
with the definition of a suitable Laplacian operator. Let
$\mathbf{T}$ be the set of transpositions of $\Omega$, the
\emph{Laplacian on $X$ associated to $\mathbf{T}$} is defined by
the operator $\Lap \in \mbox{End} (L^2(X))$ given by $$(\Lap \psi)
(x) = \sum_{g \in \mathbf{T}} \psi(g x) - |\mathbf{T}| \psi(x).$$

The Laplacian $\Lap$ belongs to the algebra $\mbox{End}_G
(L^2(X))$, this is a simple consequence of the
conjugation-invariance of $\mathbf{T}$ in $G$. Now we recall that
$\mathcal{V}_s$ denotes the representation space of $\pi_s$ for $0
\le s \le [n/2]$. Let $r_1, r_2$ and $s$ in the usual rank of
parameters and take $H: L^2(X_{r_2})_s \longrightarrow
\mathcal{V}_s$ to be an intertwining unitary isomorphism. Then we
have for $\mathrm{L}'_{\mathbf{T}} = \mathrm{L}_{\mathbf{T}} +
|\mathbf{T}| Id$
\begin{eqnarray*}
\mathrm{L}'_{\mathbf{T}} \circ \La & = & H^{-1} \circ H \circ
\mathrm{L}'_{\mathbf{T}} \circ \La = H^{-1} \circ \Big( \sum_{g
\in \mathbf{T}} \pi_s(g^{-1}) \Big) \circ H \circ \La.
\end{eqnarray*}
But the operator $\sum_{g \in \mathbf{T}} \pi_s(g^{-1})$ commutes
with the action of $\pi_s$. Therefore, by Schur's lemma, it is a
multiple of the identity $1_{\mathcal{V}_s}$. So we have shown
that the operators $\La$ are eigenvectors of the Laplacian. That
is, for each $s$ in its usual rank there exists a constant
$\mu_s(n)$ such that
\begin{equation} \label{equation11}
\Lap \circ \La + \mu_s(n) \La = 0.
\end{equation}

We recall that the classical difference operators $\delta, \Delta$
and $\nabla$ are defined by the relations $$\delta h (k) = h (k+1)
\qquad \Delta h (k) = h(k+1) - h(k) \qquad \nabla h (k) = h(k) -
h(k-1)$$ for any function $h$ of one integer variable. A
combinatorial computation shows that, if we write equation
(\ref{equation11}) in terms of the kernels, the following
difference equation is satisfied for $h_s= \lan$
\begin{eqnarray}
\nonumber & \sigma (k) \Delta \nabla h_s(k) + \tau(k) \Delta
h_s(k) + \mu_s(n) h_s(k) = 0 & \\ & \sigma (k) = k (r_1-r_2+k)
\qquad \tau(k) = r_2(n-r_1) - n k & \label{equation12} \\
\nonumber & 0 \vee (r_2-r_1) \le k \le (n-r_1) \wedge r_2 &
\end{eqnarray}

The difference equation (\ref{equation12}) is of hypergeometric
type and can be ana\-lyzed following the classical theory. With
the notation of \cite{NU}, the solutions of this equation are
given in terms of the family of \emph{Hahn polynomials}
$\widetilde{h}_s^{(\mu, \nu)} (\cdot,N)$ of parameters $N = r_1
\wedge (n-r_1) \wedge r_2 \wedge (n-r_2) + 1$, $\mu = |r_1-r_2|$,
$\nu = |n-r_1-r_2|$ and $0 \le s \le r_1 \wedge (n-r_1) \wedge r_2
\wedge (n-r_2)$.

\begin{Rem}
\emph{We observe that the left hand side of the difference
equation (\ref{equation12}) is a polynomial of degree $\le s$ for
$h_s = \lan$. Then, since the equation is satisfied at $r_1 \wedge
(n-r_1) \wedge r_2 \wedge (n-r_2) + 1$ points, it holds in fact
for all $t \in \C$.}
\end{Rem}
Now we can compute the eigenvalues $\mu_s(n)$. Namely, we already
know that equation (\ref{equation12}) can be written as $a_s t^s +
a_{s-1} t^{s-1} + \ldots + a_1 t + a_0 = 0$, hence the identity
$a_s = 0$ gives for $0 \le s \le r_1 \wedge (n-r_1) \wedge r_2
\wedge (n-r_2)$
\begin{equation} \label{equation13}
\mu_s(n) = s (n-s+1).
\end{equation}
In fact, since the numbers $\mu_s(n)$ are pairwise distinct for $0
\le s \le [n/2]$, all the eigenspaces of the linear operator
$$\begin{array}{rcl} \C^{r_1,r_2} [t] & \longrightarrow &
\C^{r_1,r_2} [t] \\ h & \longmapsto & \sigma \Delta \nabla h +
\tau \Delta h \end{array}$$ are $1$-dimensional. Here
$\C^{r_1,r_2} [t]$ stands for the space of polynomials with
complex coefficients and degree bounded by $r_1 \wedge (n-r_1)
\wedge r_2 \wedge (n-r_2)$. Therefore (\ref{equation12})
characterizes the kernel $\lan$ up to a constant factor.

The following result uses the hypergeometric equation
(\ref{equation12}) to get some useful expressions for the kernel
$\lan$. We also give the value of the highest coefficient
$a_s^{r_1,r_2}(n)$ of $\lan$. We shall use the \emph{Leibniz rule}
$\Delta (h_1 h_2) = \Delta h_1 \delta h_2 + h_1 \Delta h_2$.

\begin{The} \label{theorem3}
The kernels $\lan$ satisfy, in the usual rank of parameters, the
following relations
\begin{eqnarray*}
\lambda_s^{r_1,r_2} (t,n) & = & \sum_{j=0}^s (-1)^j
{{t}\choose{j}} \frac{[s]_j[n-s+1]_j}{[n-r_1]_j[r_2]_j} \\ \Delta
\lambda_s^{r_1,r_2} (t,n) & = & - \frac{s (n-s+1)}{r_2(n-r_1)}
\lambda_{s-1}^{r_1-1,r_2-1} (t,n-2) \\ a_s^{r_1,r_2}(n) & = &
(-1)^s \frac{[n-s+1]_s}{[n-r_1]_s[r_2]_s}.
\end{eqnarray*}
\end{The}

\begin{proof}
We start by taking $t = 0$ in equation (\ref{equation12}). Then,
by means of (\ref{equation13}) we have
\begin{equation} \label{equation14}
\Delta \lambda_s^{r_1,r_2} (0,n) = - \frac{s (n-s+1)}{r_2
(n-r_1)}.
\end{equation}
On the other hand we consider the action of the operator $\Delta$
on equation (\ref{equation12}). A calculation with the Leibniz
rule gives $$\sigma \Delta \nabla \widetilde{h}_s +
\widetilde{\tau} \Delta \widetilde{h}_s + \mu_{s-1} (n-2)
\widetilde{h}_s = 0$$ where $\widetilde{h}_s = \Delta h_s$ and
$\widetilde{\tau} (t) = (r_2-1)(n-r_1-1) - (n-2)t$. But this is
the hypergeometric equation associated to
$\lambda_{s-1}^{r_1-1,r_2-1} (\cdot,n-2)$. By (\ref{equation14})
we obtain $$\Delta \lambda_s^{r_1,r_2} (t,n) = - \frac{s
(n-s+1)}{r_2(n-r_1)} \lambda_{s-1}^{r_1-1,r_2-1} (t,n-2)$$ as we
wanted. Successive iterations of this equation give $$\Delta^j
\lambda_s^{r_1,r_2} (t,n) = (-1)^j \frac{[s]_j
[n-s+1]_j}{[n-r_1]_j[r_2]_j} \lambda_{s-j}^{r_1-j,r_2-j}
(t,n-2j)$$ for $0 \le j \le s$. Then we apply the discrete Taylor
formula $$\lambda_s^{r_1,r_2}(t,n) = \sum_{j=0}^s {{t}\choose{j}}
\Delta^j \lambda_s^{r_1,r_2}(0,n)$$ to get the desired relation.
Now, the value of $a_s^{r_1,r_2}(n)$ arises trivially.
\end{proof}

The solution of a hypergeometric equation is usually given by
Rodrigues formula, this will provide another expression for the
kernels $\lan$. The first step is to obtain the associated weight
which in the case of Hahn polynomials, up to a constant factor, is
given by
\begin{eqnarray*}
\omega^{r_1,r_2}(k,n) & = &
{{n}\choose{k,r_2-k,r_1-r_2+k,n-r_1-k}} \\ & = & \big| \big\{
(x_2,x_1) \in X_{r_2} \times X_{r_1}: |x_2 \setminus x_1| = k
\big\} \big|
\end{eqnarray*}
for $0 \vee (r_2-r_1) \le k \le r_2 \wedge (n-r_1)$. Also
$\omega^{r_1,r_2}(k,n)$ is taken to be $0$ at every integer $k$
outside that interval. The weight $\omega^{r_1,r_2}(k,n)$ is the
value at $t = k$ of the meromorphic function
$$\omega^{r_1,r_2}(t,n) = \frac{n!}{\Gamma(t+1) \Gamma(r_2-t+1)
\Gamma(r_1-r_2+t+1) \Gamma(n-r_1-t+1)}$$ which satisfies the
relations
\begin{eqnarray}
\Delta (\omega^{r_1,r_2} (t,n) \sigma (t)) & = & \omega^{r_1,r_2}
(t,n) \tau (t) \label{equation15} \\ \delta (\omega^{r_1,r_2}
(t,n) \sigma (t)) & = & n (n-1) \omega^{r_1-1,r_2-1} (t,n-2)
\label{equation16}
\end{eqnarray}
Leibniz rule and (\ref{equation15}) give the self-adjoint form of
the hypergeometric equation (\ref{equation12}) $$\Delta
(\omega^{r_1,r_2} (t,n) \sigma (t) \nabla \la (t,n)) + \mu_s(n)
\omega^{r_1,r_2} (t,n) \la (t,n) = 0.$$ On the other hand we can
combine Leibniz rule, the recurrence given in theorem
\ref{theorem3} and (\ref{equation16}) to give, from the
self-adjoint form of (\ref{equation12}), the relation
$$\omega^{r_1,r_2} (t,n) \la (t,n) = \frac{n (n-1)}{(n-r_1) r_2}
\nabla (\omega^{r_1-1,r_2-1} \lambda_{s-1}^{r_1-1,r_2-1})
(t,n-2).$$ By induction we get, for $0 \le s \le r_1 \wedge
(n-r_1) \wedge r_2 \wedge (n-r_2)$ and $t \in \C$, the Rodrigues
formula
\begin{equation} \label{equation17}
\omega^{r_1,r_2}(t,n) \la (t,n) =
\frac{[n]_{2s}}{[n-r_1]_s[r_2]_s} \nabla^s
\omega^{r_1-s,r_2-s}(t,n-2s).
\end{equation}

\noindent Rodrigues formula provides some alternative expressions
for the kernels $\lan$.
\begin{enumerate}
\item Taking into account Newton's binomial formula
$$\nabla^s = (1 - \delta^{-1})^s = \sum_{j=0}^s (-1)^j
{{s}\choose{j}} \delta^{-j}$$ we easily obtain from
(\ref{equation17}) the relation $$\lambda_s^{r_1,r_2} (t,n) =
\sum_{j=0}^s \frac{(-1)^j}{[n-r_1]_s[r_2]_s} {{s}\choose{j}}
[t]_j[r_2-t]_{s-j}[r_1-r_2+t]_j[n-r_1-t]_{s-j}$$ for $r_1, r_2$
and $s$ in the usual rank of parameters and $t \in \C$.
\item As a particular case, if we evaluate this relation at
the integer variable $k$, we get the classical expressions
\begin{eqnarray*} \lambda_s^{r_1,r_2}(k,n) & = &
\frac{[r_1]_s}{[r_2]_s} \sum_{j=0}^s (-1)^j \frac{{{s}\choose{j}}
{{r_1-s}\choose{r_1-r_2+k-j}}
{{n-r_1-s}\choose{k-j}}}{{{r_1}\choose{r_1-r_2+k}}
{{n-r_1}\choose{k}}} \\ \lambda_s^{r_1,r_2}(k,n) & = &
\frac{[n-r_2]_s}{[n-r_1]_s} \sum_{j=0}^s (-1)^j
\frac{{{s}\choose{j}} {{r_2-s}\choose{k-j}}
{{n-r_2-s}\choose{r_1-r_2+k-j}}}{{{r_2}\choose{k}}
{{n-r_2}\choose{r_1-r_2+k}}}.
\end{eqnarray*}
\end{enumerate}

\section{The Fourier transform $\mathcal{F}_X$}
\label{section5}

\setcounter{The}{0}

Given $x, x' \in X$ we can consider the symmetric difference
operator $x \triangle x' = (x \cup x') \setminus (x \cap x')$
which provides a structure of abelian group on $X$. The dual group
of $X$ is given by the set $\widehat{X} = \{\chi_x\}_{x \in X}$ of
characters of $X$, where $\chi_x (x') = (-1)^{|x \cap x'|}$. This
group structure on $X$ leads us to consider the Fourier transform
$\mathcal{F}_X: L^2(X) \longrightarrow L^2(\widehat{X})$, which is
given by the formula $$(\mathcal{F}_X \psi) (\chi_x) = \sum_{x'
\in X} (-1)^{|x \cap x'|} \psi(x').$$ Moreover, we can see the
Fourier transform $\mathcal{F}_X$ as an element of the algebra
$\mbox{End} (L^2(X))$. In fact it is not difficult to check that
$\mathcal{F}_X \in \mbox{End}_G (L^2(X))$. Therefore
$\mathcal{F}_X$ can be written as a linear combination of our
basis. This time we shall make use of the unitary operators
$\widetilde{\Lambda}_s^{r_1,r_2}$, this will simplify some of the
results in the sequel. Before investigating the coefficients of
this linear combination, we introduce a family of orthogonal
polynomials of hypergeometric type. The \emph{Krawtchouk
polynomials} $P_m (k,n)$ ($0 \le m \le n$) are defined as the
solutions of the hypergeometric equation $k \nabla \Delta h (k) +
(n-2k) \Delta h (k) + \mu_m h (k) = 0$ for $0 \le k \le n$,
normalized by the condition $P_m (0,n) = {{n}\choose{m}}$. The
Krawtchouk polynomials have the form $$P_m (k,n) = \sum_{j=0}^m
(-1)^j {{k}\choose{j}} {{n-k}\choose{m-j}}$$ and satisfy the
following relation for $0 \le k \le n$ and $0 \le m \le n$ $$P_m
(k,n) = (-1)^m P_m (n-k,n).$$ A more detailed exposition of these
topics can be found in \cite{MS} and \cite{NU}.

\begin{The} \label{theorem4}
The Fourier transform on $X$ satisfies the following decomposition
$$\mathcal{F}_X = \sum_{0 \le r_1,r_2 \le n}
\sum_{s=0}^{N(r_1,r_2)} k_s^{r_1,r_2}
\widetilde{\Lambda}_s^{r_1,r_2}$$ where $N(r_1,r_2) = r_1 \wedge
(n-r_1) \wedge r_2 \wedge (n-r_2)$ and $$k_s^{r_1,r_2} = (-2)^s
\sqrt{{{n-2s}\choose{r_2-s}} \Big/ {{n-2s}\choose{r_1-s}}}
P_{r_1-s} (r_2-s,n-2s).$$
\end{The}

\begin{proof} Using the pairwise orthogonality of the operators
$\widetilde{\Lambda}_s^{r_1,r_2}$ with respect to the
Hilbert-Schmidt product we obtain $$k_s^{r_1,r_2} =
\frac{\mbox{tr} \big( \mathcal{F}_X \circ
(\widetilde{\Lambda}_s^{r_1,r_2})^{\star} \big)}{\mbox{tr} \big(
\widetilde{\Lambda}_s^{r_1,r_2} \circ
(\widetilde{\Lambda}_s^{r_1,r_2})^{\star} \big)} =
\frac{1}{{{n}\choose{s}} - {{n}\choose{s-1}}} \ \ \mbox{tr}
(\mathcal{F}_X \circ \widetilde{\Lambda}_s^{r_2,r_1}).$$ where we
have used the relations for the operators
$\widetilde{\Lambda}_s^{r_1,r_2}$ of section \ref{section3}. But,
recalling that $\omega^{r_1,r_2}(k,n)$ is the cardinal of the
orbit $\theta_k$ for $k$ in the usual interval and zero otherwise,
we can write
\begin{eqnarray*}
\mbox{tr} (\mathcal{F}_X \circ \widetilde{\Lambda}_s^{r_2,r_1}) &
= & \alpha_s^{r_2,r_1} \sum_{x_1 \in X_{r_1}} \sum_{x_2 \in
X_{r_2}} (-1)^{|x_1 \cap x_2|} \lambda_s^{r_2,r_1} (|x_1 \setminus
x_2|,n) \\ & = & \alpha_s^{r_2,r_1} \sum_{t \in \Z} (-1)^{r_1 - t}
\omega^{r_2,r_1}(t,n) \lambda_s^{r_2,r_1}(t,n)
\end{eqnarray*}
where we recall that $\alpha_s^{r_1,r_2}$ denotes the
normalization constant of (\ref{equation10}). Then we apply
Rodrigues formula (\ref{equation17}) and summation by parts to get
\begin{eqnarray*}
\mbox{tr} (\mathcal{F}_X \circ \widetilde{\Lambda}_s^{r_2,r_1}) &
= & \alpha_s^{r_2,r_1} \frac{(-1)^{r_1}
[n]_{2s}}{[r_1]_s[n-r_2]_s} \sum_{t \in \Z} (-1)^t \nabla^s
\omega^{r_2-s,r_1-s} (t,n-2s) \\ & = & \alpha_s^{r_2,r_1}
\frac{(-1)^{r_1} 2^s [n]_{2s}}{[r_1]_s[n-r_2]_s} \sum_{t \in \Z}
(-1)^t \omega^{r_2-s,r_1-s} (t,n-2s).
\end{eqnarray*}
Now, the relation given in section \ref{section4} for the weight
$\omega^{r_1,r_2} (\cdot,n)$ in terms of binomial numbers gives
the following expression for $\mbox{tr} (\mathcal{F}_X \circ
\widetilde{\Lambda}_s^{r_2,r_1})$ $$\frac{(-1)^{r_1} 2^s [n]_{2s}
{{n-2s}\choose{r_2-s}}}{\sqrt{{{n-2s}\choose{r_1-s}}
{{n-2s}\choose{r_2-s}}} s! [n-s+1]_s} \sum_{t \in \Z} (-1)^t
{{r_2-s}\choose{r_2-r_1+t}} {{n-r_2-s}\choose{t}}.$$ Then we use
the properties of the Krawtchouk polynomials to obtain $$\mbox{tr}
(\mathcal{F}_X \circ \widetilde{\Lambda}_s^{r_2,r_1}) =
\frac{(-2)^s [n]_{2s}}{s! [n-s+1]_s} \sqrt{{{n-2s}\choose{r_2-s}}
\Big/ {{n-2s}\choose{r_1-s}}} P_{r_1-s} (r_2-s,n-2s).$$ Dividing
by ${{n}\choose{s}} - {{n}\choose{s-1}}$ we obtain the desired
relation.
\end{proof}

\begin{Rem}
\emph{The self-adjointness of the Fourier transform
$\mathcal{F}_X$ is now a simple consequence of theorem
\ref{theorem4}. The identity $k_s^{r_1,r_2} = k_s^{r_2,r_1}$
arises as one of the relations that characterize the Krawtchouk
polynomials, namely $${{n}\choose{k}} P_m (k,n) = {{n}\choose{m}}
P_k (m,n).$$ Thus we just need to express $\mathcal{F}_X$ in terms
of the operators $\widetilde{\Lambda}_s^{r_1,r_2}$ as in theorem
\ref{theorem4} and apply the relation
$(\widetilde{\Lambda}_s^{r_1,r_2})^{\star} =
\widetilde{\Lambda}_s^{r_2,r_1}$ to obtain the identity
$\mathcal{F}_X = \mathcal{F}_X^{\star}$.}
\end{Rem}
We close this section with an expression of the Fourier transform
$\mathcal{F}_X$ which reflects the action of the representation
$\rho$. For that purpose we start by decomposing $L^2(X)$ into
irreducible subspaces as follows $$L^2(X) = \bigoplus_{r=0}^n
\bigoplus_{s=0}^{r \wedge (n-r)} L^2(X_r)_s =
\bigoplus_{s=0}^{[n/2]} \bigoplus_{r=s}^{n-s} L^2(X_r)_s.$$ But we
know that the spaces $L^2(X_r)_s$ are $G$-equivalent to
$\mathcal{V}_s$, the representation space of $\pi_s \in \G$, for
$s \le r \le n-s$. Therefore, if we declare $\mathcal{V}_s =
L^2(X_s)_s$ we have the following unitary intertwining isomorphism
$$\begin{array}{crcl} T: & \displaystyle \bigoplus_{s=0}^{[n/2]}
\mathcal{V}_s^{n-2s+1} & \longrightarrow & L^2(X) \\ &
(\mathrm{v}_s^r) & \longmapsto & \displaystyle \sum_{s=0}^{[n/2]}
\sum_{r=s}^{n-s} \widetilde{\Lambda}_s^{s,r} (\mathrm{v}_s^r).
\end{array}$$ Now we define $\mathrm{F}_X = T^{-1} \circ
\mathcal{F}_X \circ T$. It is obvious that $\mathrm{F}_X$ will
give an expression of $\mathcal{F}_X$ very related to the action
of the representation $\rho$. The following theorem explains that
relation.

\begin{The} \label{theorem5}
The operator $\mathrm{F}_X$ has the form $\mathrm{F}_X
(\mathrm{v}_s^r) = (\mathrm{w}_s^r)$ where
$$\mathrm{w}_{s_0}^{r_0} = \sum_{r=s_0}^{n-s_0} k_{s_0}^{r,r_0}
\mathrm{v}_{s_0}^r$$ and the numbers $k_{s_0}^{r,r_0}$ are the
coefficients of $\mathcal{F}_X$ given in theorem $\ref{theorem4}$.
\end{The}

\begin{proof}
Let us observe that $$\sum_{s=0}^{[n/2]} \sum_{r=s}^{n-s}
\widetilde{\Lambda}_s^{s,r} (\mathrm{w}_s^r) = T \circ
\mathrm{F}_X (\mathrm{v}_s^r) = \mathcal{F}_X \circ T
(\mathrm{v}_s^r) = \sum_{s=0}^{[n/2]} \sum_{r=s}^{n-s}
\mathcal{F}_X \circ \widetilde{\Lambda}_s^{s,r}
(\mathrm{v}_s^r).$$ Hence, by projecting onto the isotipic
component correspondent to $s_0$, we can fix the value of $s$
obtaining $$\sum_{r=s_0}^{n-s_0} \widetilde{\Lambda}_{s_0}^{s_0,r}
(\mathrm{w}_{s_0}^r) = \sum_{r={s_0}}^{n-s_0} \mathcal{F}_X \circ
\widetilde{\Lambda}_{s_0}^{s_0,r} (\mathrm{v}_{s_0}^r).$$ Let us
fix $r = r_0$, then left multiplication by
$\widetilde{\Lambda}_{s_0}^{r_0,s_0}$ and theorem \ref{theorem4}
imply $$\mathrm{w}_{s_0}^{r_0} = \sum_{r=s_0}^{n-s_0}
\widetilde{\Lambda}_{s_0}^{r_0,s_0} \circ \mathcal{F}_X \circ
\widetilde{\Lambda}_{s_0}^{s_0,r} (\mathrm{v}_{s_0}^r) =
\sum_{r=s_0}^{n-s_0} k_{s_0}^{r,r_0} \mathrm{v}_{s_0}^r$$ as we
wanted. This completes the proof.
\end{proof}

\begin{Rem} \emph{If $\mathbf{M}(k)$ stands for the algebra of
$k \times k$ matrices with complex entries, then we can see the
space $\mbox{End}_G (L^2(X))$ as a matrix algebra via the
following algebra isomorphism $$\begin{array}{crcl} \Phi: &
\displaystyle \bigoplus_{s=0}^{[n/2]} \mathbf{M}(n-2s+1) &
\longrightarrow & \mbox{End}_G (L^2(X)) \\ & \Big( \,\ a_s \,\
\Big)_{0 \le s \le [n/2]} & \longmapsto & \displaystyle
\sum_{s=0}^{[n/2]} \sum_{s \le r_1,r_2 \le n-s} a_s^{r_1,r_2}
\widetilde{\Lambda}_s^{r_1,r_2}
\end{array}$$ where $a_s = \Big( \,\ a_s^{r_1,r_2} \Big)$ with $s
\le r_1, r_2 \le n-s$. Now we define $$\mathcal{K} = \Big( \,\ k_s
\,\ \Big) \in \bigoplus_{s=0}^{[n/2]} \mathbf{M}(n-2s+1) \qquad
\mbox{such that} \qquad k_s = \Big( \,\ k_s^{r_1,r_2} \,\ \Big)$$
where the numbers $k_s^{r_1,r_2}$ denote the coefficients of
$\mathcal{F}_X$ given in theorem \ref{theorem4}. Then we have that
\begin{itemize}
\item[$(a)$] Theorem \ref{theorem4} asserts that $\mathcal{F}_X =
\Phi (\mathcal{K})$.
\item[$(b)$] Theorem \ref{theorem5} asserts that the matrix of
$\mathrm{F}_X$ is given by $\mathcal{K}^{\mbox{t}}$.
\end{itemize}}
\end{Rem}

\bibliographystyle{amsplain}

\end{document}